\documentclass[12pt, reqno]{amsart}
\usepackage{amsfonts}
\usepackage{amsmath}
\usepackage{amssymb}
\usepackage{eufrak}
\usepackage{mathrsfs}
\usepackage{bbding}
\usepackage{fancyhdr}

\numberwithin{equation}{section} 

\input xy
\xyoption{all}
\title{On algebraic equivalences among the 27 Abel-Prym curves on a generic abelian 5-fold}
\author{Maxim Arap}
\author{Robert Varley}

\begin{document}

\begin{abstract}
This article shows that on a generic  principally polarized abelian variety of dimension five the $\mathbb{Q}$-vector space of algebraic equivalences among the 27 Abel-Prym curves has dimension $ 20$. 
\end{abstract}

\address{Johns Hopkins University, Department of Mathematics,  404 Krieger Hall, 3400 N. Charles Street, Baltimore, MD 21218, USA}
\email{marap@math.jhu.edu}

\address{
University of Georgia,
Department of Mathematics,
Boyd Graduate Studies Research Center,
Athens, GA 30602,
USA}
\email{rvarley@math.uga.edu}

\newcommand{\bP}{\mathbb{P}}
\newcommand{\bC}{\mathbb{C}}
\newcommand{\mc}{\mathcal}
\newcommand{\ra}{\rightarrow}
\newcommand{\thra}{\twoheadrightarrow}
\newcommand{\mb}{\mathbb}
\newcommand{\mrm}{\mathrm}
\newcommand{\p}{\prime}
\newcommand{\ms}{\mathscr}
\newcommand{\pl}{\partial}
\newcommand{\ti}{\tilde}
\newcommand{\wti}{\widetilde}
\newcommand{\ol}{\overline}
\newcommand{\ul}{\underline}
\newcommand{\Sp}{\mrm{Spec}\,}

\newtheorem{lem1}{Lemma}[section]
\newtheorem{lem2}[lem1]{Lemma}
\newtheorem{lem3}[lem1]{Lemma}
\newtheorem{lem4}[lem1]{Lemma}
\newtheorem{lem5}[lem1]{Lemma}
\newtheorem{lem6}[lem1]{Lemma}
\newtheorem{def1}[lem1]{Definition}
\newtheorem{thm1}[lem1]{Theorem}
\newtheorem*{thm0}{Theorem}
\newtheorem{pro1}[lem1]{Proposition}
\newtheorem{rem1}[lem1]{Remark}
\newtheorem{rem2}[lem1]{Remark}
\newtheorem{cor1}[lem1]{Corollary}

\renewcommand{\thefootnote}{}

\maketitle
 \footnote{\today}

\thispagestyle{empty}

\section{Introduction}
Let $X$ be a principally polarized abelian variety (ppav) over $\mb{C}$. Let $A^*(X)$ denote the Chow ring of $X$ modulo \emph{algebraic} equivalence with $\mb{Q}$-coefficients (this notation differs from the one in  \cite[10.3, p.185]{Fu98}). Besides the intersection product, the ring $A^*(X)$ is endowed with Pontryagin product defined by $$x_1 \ast x_2 = m_\ast (p_1^\ast x_1 \cdot p_2^\ast x_2),$$
where $m\colon X \times X \ra X$ is the addition morphism, and $p_j\colon X\times X \ra X$ is the projection onto the $j^{\mrm{th}}$ factor, cf. [BL, p.530]. Moreover, $A^*(X)$ carries a bi-grading,$$A(X) = \bigoplus_{l,s}A^l(X)_{(s)}.$$  The $l$-grading is by codimension. The Beauville grading ($s$) is defined by the condition $x \in A^l(X)_{(s)}$ if and only if $k^*x=k^{2l-s}x$ for all $k \in \mb{Z}$, where $k$ also denotes the endomorphism of $X$ given by $x \mapsto kx$ (the second grading was originally defined in \cite{B86} for rational equivalence but it is well-defined for algebraic equivalence as well). The $(s)$-component of a cycle $Z$ is denoted by $Z_{(s)}$. Also, let us recall the Fourier transform $\mc{F}_X\colon A^*(X) \ra A^*(X)$ that is given by $z\mapsto p_{2,*}(p_1^*z\cdot e^\ell)$, where $\ell$ is the class of the Poincar\'e bundle on $X\times X$ (here we identify $X$ and its dual abelian variety using the principal polarization). 

Let  $P:=\mrm{Prym}(\ti{C}/C)$ be the Prym variety of a connected \'etale double cover $\ti{C} \ra C$ of smooth curves, see \cite{Mu74}. In the sequel we assume that $\ti{C}$ is not hyperelliptic and we let $p:=\dim P$. After fixing a base-point $\ti{o} \in \ti{C}$ we have the Abel-Prym map $\psi\colon  \ti{C} \ra P$ whose image will be denoted by $\ti{C}$ as well and is called an Abel-Prym curve ($\psi$ is a closed embedding, cf. \cite[Cor.12.5.6, p.380]{BL04}). 

Given an Abel-Prym curve $\ti{C} \subset P$, by \cite[Def.1, p.708]{A12} there is an associated tautological (sub)ring $\ms{T}(P, \ti{C}) \subset A^*(P)$, defined to be the smallest subring of $A^*(P)$ for the intersection product that contains $[\ti{C}]$ and is stable under $*, \mc{F}_P$ and $k^*$ for all $k \in \mb{Z}$. By \cite[Thm.4, p.710]{A12}, $\ms{T}(P, \ti{C})$ is generated by the cycles $\zeta_n:=\mc{F}_P([\ti{C}]_{(n-1)})$ for $1 \le n \le p-1$ odd (the components of $[\ti{C}]$ of even Beauville degree vanish because $\ti{C}$ has a symmetric translate). Also, by \cite[Rem.3, p.711]{A12} the ring $\ms{T}(P, \ti{C})$ does not depend on the choice of an Abel-Prym curve if $P$ is a generic ppav of dimension $p \ne 5$. As a corollary of the main result of this article (Theorem \ref{thm1}), we obtain that for a generic ppav $P$ of dimension $5$ the tautological ring $\ms{T}(P, \ti{C})$ does depend on the choice of an Abel-Prym curve $\ti{C} \subset P$. This gives $27$ tautological rings,  one for each choice of an Abel-Prym curve (see the next paragraph for the explanation of the number $27$). However, by \cite[p.35]{A11} each of these $27$ rings is isomorphic to the quotient $\mb{Q}[x_1,x_3]/(x_1^6 , x_3^2 , x_1^2x_3)$ of the polynomial ring in two variables via the homomorphism $\zeta_i \mapsto x_i$ (the cycles $\zeta_n$ vanish for $n > 3$ by \cite[Cor.3.3.8, p.34]{A11}).

Let $\ms{P}\colon \ms{R}_6 \ra \ms{A}_5$ be the Prym map from the moduli space $\ms{R}_6$ of connected \'etale double covers $\ti{C}\ra C$ of smooth genus 6 curves to the moduli space $\ms{A}_5$ of 5-dimensional ppav's. By \cite{DS81}, the morphism $\ms{P}$ is generically finite of degree $27$. By the tetragonal construction described in \cite[2.5, p.76]{Do92}, given $(\ti{C}\ra C) \in \ms{R}_6$ and a degree 4 linear pencil  $\frak{g}^1_4$ on $C$ there exist \'etale double covers $\ti{C}_0/C_0$ and  $\ti{C}_1/C_1$ of tetragonal curves whose associated Prym is isomorphic (as a ppav) to $P=\mrm{Prym}(\ti{C}/C)$. Furthermore, by \cite[Thm.4.2, p.90]{Do92}, the correspondence on the generic fiber of $\ms{P}$ induced by the tetragonal construction is isomorphic to the incidence correspondence of the $27$ lines on a smooth cubic surface. Also, the monodromy group of the cover $\ms{R}_6\ra \ms{A}_5$ is the Weyl group $W$ of the $E_6$ lattice. 

It follows from \cite{Fa96} and \cite{A12} that the $27$ Abel-Prym curves on a generic ppav of dimension $5$ are not pairwise algebraically equivalent, see Lemma \ref{lem1}. In this article we determine all the algebraic equivalences (with $\mb{Q}$-coefficients) among these $27$ effective 1-cycles. The main result is the following. 

\begin{thm0}[=Theorem \ref{thm1}] \label{thm0}
The $\mb{Q}$-vector space of algebraic equivalences of the $27$ Abel-Prym curves on a generic principally polarized abelian 5-fold has dimension $20$.
\end{thm0}

The key ingredient for proving the above theorem is Lemma \ref{lem4}. This lemma implies that the action of the monodromy group of the cover $\ms{R}_6 \ra \ms{A}_5$ on the $27$ Abel-Prym curves on a generic ppav of dimension $5$ preserves algebraic equivalences among the Abel-Prym curves  (see \cite{Har79} for a discussion of monodromy in algebraic geometry and its applications).

Throughout the paper we work over the field of complex numbers. 

\section{Algebraic equivalences}

\begin{pro1}[Connectedness principle] \label{pro1}
Let $Y$ be a smooth connected but not necessarily complete curve and let $f\colon X \ra Y$ be a proper morphism. If $X$ is connected then the locus 
$$\mrm{CF}(f):=\{y \in Y \, | \, X_y \text{ is connected} \}$$
parametrizing connected fibers is Zariski closed in $Y$. 
\end{pro1}
\begin{proof}  Let $f=h \circ g$ be the Stein factorization, where $g \colon X \ra Z$ has connected fibers and $h\colon Z \ra Y$ is a finite morphism. Since $Z=g(X)$ and $X$ is connected then $Z$ is connected. If $Z$ is a point, the conclusion holds trivially. Assume that $Z$ is not a point, let $Z_1,\ldots, Z_k$ be the irreducible components of $Z$ and let $h_i$ be the restriction of $h$ to $Z_i$. If $k \ge 2$ then there are at most finitely many points with connected fibers since each $Z_i$ is mapping onto $Y$.  Thus it remains to consider the case $k =1$, i.e., the case when  $Z$ is an irreducible curve and $h\colon Z \twoheadrightarrow Y$ is a finite surjective morphism. We may and shall assume that $Z$ is reduced. Consider the normalization $\tilde{Z} \to Z$, which is a finite birational morphism onto $Z$, and the composed morphism $\ti{h}\colon \tilde{Z} \to Z \to Y$.  It suffices to show that the locus $\mrm{CF}(\ti{h})$ is Zariski closed in $Y$.  If $\ti{h}$ has degree $1$ then $\ti{h}$ is an isomorphism and $\mrm{CF}(\ti{h}) = Y$.  If $\ti{h}$ has degree $\ge 2$, then $\mrm{CF}(\ti{h})$, which must be contained in the branch locus, is at most a finite set of points and hence is Zariski closed in $Y$.
\end{proof}

\begin{rem1}
The conclusion of the above proposition becomes false in general if  $Y$ is a reducible curve, $\dim Y > 1$ or if $X$ is not connected.
\end{rem1}

\begin{lem2}\label{lem2}
If $\ti{C}_1/C_1,\ti{C}_2/C_2, \ti{C}_3/C_3$ is a tetragonally related triple on a Prym variety $P$ then $[\ti{C}_1]_{(2)}+[\ti{C}_2]_{(2)}+[\ti{C}_3]_{(2)}=0$ in $A^*(P)$. 
\end{lem2}
\begin{proof} Taking $\ti{C}_1/C_1$ as a reference double cover, the curves $\ti{C}_2, \ti{C}_3 \subset P$ are the special subvarieties associated to some $\frak{g}^1_4$ on $C_1$, see \cite[3(c), p.366]{B82}. By \cite[Ex.1, p.724]{A12} we have the relation $[\ti{C}_2]_{(2)}+[\ti{C}_3]_{(2)}= -[\ti{C}_1]_{(2)}$.
\end{proof}

\begin{lem1} \label{lem1}
On a generic ppav of dimension 5 there is a pair of Abel-Prym curves that are not algebraically equivalent. 
\end{lem1}
\begin{proof}
By \cite{Fa96}, for a generic ppav $P \in \ms{A}_5$ with Prym realization $P=\mrm{Prym}(\ti{C}/C)$ we have $[\ti{C}]_{(2)}\ne 0$. It is well known that $C$ has five $\frak{g}^1_4$'s. Choosing a $\frak{g}^1_4$ on the base curve $C$,  the special subvarieties $V_0, V_1 \subset P$ associated to the $\frak{g}^1_4$ are also Abel-Prym curves, \cite[3(c), p.366]{B82}.  By Lemma \ref{lem2} we have the relation $[V_0]_{(2)}+[V_1]_{(2)} + [\ti{C}]_{(2)}=0$. If the Abel-Prym curves $V_0,V_1$ and $\ti{C}$ are pairwise algebraically equivalent then we have $3[\ti{C}]_{(2)}=0$. Since we are working with $\mb{Q}$-coefficients  this implies $[\ti{C}]_{(2)}=0$, which is a contradiction. 
\end{proof}

\begin{lem3} \label{lem3}
If $\ti{C}_1,  \ldots, \ti{C}_{27}$ are the $27$ Abel-Prym curves on a generic ppav $P\in \ms{A}_5$ then the vector space of algebraic equivalences among $\ti{C}_1,  \ldots, \ti{C}_{27}$ induced by the tetragonal construction has dimension 20. 
\end{lem3}
\begin{proof}
Let $\mc{V}$ be the $\mb{Q}$-vector space with basis $\ti{C}_1, \ldots, \ti{C}_{27}$ and consider the linear map $$cl \colon \mc{V} \ra A^4(P)$$ that takes an element of $\mc{V}$ to its class modulo algebraic equivalence. Let $\mc{V}_{tet} \subset \mc{V}$ be the subspace spanned  by the vectors $\ti{C}_i+\ti{C}_j+\ti{C}_k$, where $\ti{C}_i,\ti{C}_j,\ti{C}_k$ are tetragonally related. The vector space $R_{tet}:=\ker(cl \colon \mc{V}_{tet} \ra A^4(P))$ is the space of algebraic equivalences induced by the tetragonal construction.

It follows immediately from \cite[Thm.4.2, p.90]{Do92} that tetragonally related triples are in bijection with the 45 triples of lines forming triangles on a smooth cubic surface in $\bP^3$. Let us fix an ordering $L_1, \ldots, L_{27}$ of the 27 lines on the cubic surface. Each of the triangles determines a vector in a 27 dimensional $\mb{Q}$-vector space with basis $L_1, \ldots, L_{27}$ with 1's in the coordinates corresponding to the lines forming the triangle, see \cite{Dol10} and \cite[Ch.4]{Hun96} for details on the geometry of cubic surfaces. 
An explicit calculation shows that the dimension of the span of these 45 vectors is 21. Therefore, $\dim \mc{V}_{tet} = 21$. 

For each $i \in\{1, \ldots, 27\}$ we have the formula $[\ti{C}_i] = [\ti{C}_i]_{(0)}+[\ti{C}_i]_{(2)}$ in $A^4(P)$.  The vanishing of the odd degree Beauville components of $[\ti{C}_i]$ is a consequence of the fact that $\ti{C}_i$ has a symmetric translate. The vanishing of components of degree $\ge 4$ was proven in \cite[Cor.3.3.8]{A11}. Furthermore, using the Fourier transform on $A^*(P)$ it is easy to see that $[\ti{C}_1]_{(0)}=\ldots =[\ti{C}_{27}]_{(0)}$ in $A^4(P)$. The image of $[\ti{C}_i]_{(0)}$ under the cycle class map $A^4(P) \ra H^8(P, \mb{Q})$ is twice the minimal cohomology class, and therefore, $[\ti{C}_i]_{(0)} \ne 0$ for all $i\in \{1,\ldots,27\}$. Using Lemma \ref{lem2}, this implies that for every tetragonally related triple $\ti{C}_i,\ti{C}_j,\ti{C}_k$ the image vector $cl(\ti{C}_i+\ti{C}_j+\ti{C}_k)$ lies in the 1-dimensional subspace of $A^4(P)$ spanned  by $[C_1]_{(0)}$. Therefore, $\dim R_{tet}=20$.  
\end{proof} 

Let $L_1, \ldots, L_{27}$ be the lines on a smooth cubic surface. It is well known that the group of symmetries of these lines is the Weyl group $W$ of the $E_6$ lattice. 
\begin{lem5}\label{lem5}
The representation of $W$ on $\bigoplus_{i=1}^{27}\mb{Q}L_i$ decomposes into a direct sum of irreducible representations of dimensions 1,6 and 20. 
\end{lem5}
\begin{proof} The Weyl group $W$ has order $51840$ and $25$ conjugacy classes. Computing the character table for $W$ and using the orthogonality relations between irreducible characters, it can be seen that $\bigoplus_{i=1}^{27}\mb{Q}L_i$ decomposes into irreducible representations of dimensions $1,6$ and $20$. Our calculations were carried out on {\tt Magma}.
\end{proof}

\begin{rem2}
The vectors determined by the $45$ triangles in the cubic surface span a $21$-dimensional subrepresentation, whose quotient is the representation of $W$ on the $\mb{Q}$-span of the $E_6$ lattice, which is well-known to be absolutely irreducible, \cite[Lem.B, p.53]{Hum73}. Under the map of the free abelian group on $L_1, \ldots, L_{27}$ to the Picard group of the cubic surface, each triangle is mapped to $-K$ (the anti-canonical divisor class) and the perp of $K$ is precisely the $E_6$ lattice. 
\end{rem2}

Let $\ms{A}_5^\circ$ be the complement in $\ms{A}_5$ of the union of the branch divisor of $\ms{R}_6 \ra \ms{A}_5$ and the singular locus of $\ms{A}_5$. Let $V$ be a non-empty connected Zariski open subset of $\ms{A}_5^\circ$ over which the Prym map $\ms{R}_6 \ra \ms{A}_5$ is finite. We have the following key lemma. 

\begin{lem6} \label{lem6}
A general (smooth, connected but not complete) curve section $S$ of $V$  has the property that the map $\pi_1(S) \ra \pi_1(V)$ induced by the inclusion is surjective. 
\end{lem6}
\begin{proof}
If $V$ were projective, the result of the lemma would follow immediately by the Lefschetz hyperplane theorem for homotopy groups. Since $V$ is  not projective, we have to compactify $V$ and resolve the resulting space before applying Lefschetz theorem. Although the argument is standard, we include it here for completeness.

Let $\bar{V}$ be a smooth projective variety containing $V$ and such that the complement $\bar{V} - V$ is a simple normal crossings (snc) divisor. To obtain $\bar{V}$, we may take the closure of $V$ in some projective space and apply Hironaka's theorem to resolve the singularities and make the complement of $V$ an snc divisor (we identify $V$ with its proper transform in  $\bar{V}$).    

Let $D_1, \ldots, D_n$ be the irreducible components of $\bar{V} - V$. 
For a general curve section $\bar{S}$ of $\bar{V}$ we have $\bar{S} \cap D_i \ne \emptyset$ for $1\le i \le n$  and $S:=\bar{S} \cap V$ is a non-empty Zariski open set in $\bar{S}$. Also, by the theorem of Lefschetz (cf. \cite[Thm.7.4, p.41]{Mi63}) and induction we have a surjection $\pi_1(\bar{S}) \twoheadrightarrow \pi_1(\bar{V})$. 

Given a smooth connected complex manifold $M$ and a smooth connected divisor $\Gamma \subset M$, using van Kampen's theorem we may check that there is an exact sequence $$\mb{Z}\ell_\Gamma \ra \pi_1(M - \Gamma) \ra \pi_1(M) \ra 1,$$
where $\ell_\Gamma$ is the class of a loop in $M$ that goes once around $\Gamma$. 

Since $\bar{S} \cap D_1 \ne \emptyset$ then $\ell_{D_1}$ is in the image of $\pi_1(S) \ra \pi_1(\bar{V}- D_1)$. Therefore, using the above exact sequence we have a natural diagram with all maps being surjective 
$$\xymatrix{
\pi_1(S) \ar@{->>}[r] \ar@{->>}[d] & \pi_1(\bar{S}) \ar@{->>}[d] \\
\pi_1(\bar{V}- D_1) \ar@{->>}[r] & \pi_1(\bar{V}).
}$$ 
Now, removing $D_i$'s one at a time and applying the above argument, we conclude by induction that there is a surjection $\pi_1(S) \twoheadrightarrow \pi_1(V)$. 
\end{proof}

Recall that a connected cover $f \colon X \ra Y$ of topological spaces is said to be Galois (or regular, or normal, see \cite[p.70]{Hat02}) if for every $y \in Y$ and every pair of lifts $x_1, x_2 \in X$ of $y$ there is a deck transformation of $X$ taking $x_1$ to $x_2$.  Given a connected cover $f \colon X \ra Y$ of topological spaces we shall denote its Galois closure by $\ti{X} \ra Y$. More precisely, $\ti{X}$ is the connected cover of $Y$ corresponding to the normal subgroup of $\pi_1(Y)$ obtained as the intersection of all the conjugates of $f_*(\pi_1(X))$ in $\pi_1(Y)$. Also,  the group of deck transformations of $\ti{X}$ over $Y$ will be referred to as the monodromy group of $X$ over $Y$ and will be denoted by $M(X/Y)$ in the sequel. 

Let $U:=\ms{P}^{-1}(V)$ be the inverse image of $V$ under the Prym map $\ms{P}\colon \ms{R}_6 \ra \ms{A}_5$. By our construction, the cover $U \ra V$ is finite \'etale of degree 27. Unfortunately $V$ does not carry a universal family of ppav's, which we shall need in the sequel. In what follows, we shall show that there is a finite cover $V^\p \ra V$ such that $V^\p$ carries a universal family of ppav's and the monodromy group $M(U^\p/V^\p)$ of the pull-back $U^\p \ra V^\p$ of $U\ra V$ surjects onto $W$. Let us argue that taking ppav's with level 3 structure suffices for this purpose, i.e., we may take the connected \'etale cover $\lambda \colon V^\p \ra V$ with Galois group $G:=\mrm{Sp}(10, \mb{Z}/3)/(\pm I)$. First, it is well known that the moduli space of 5-dimensional ppav's with level 3 structure carries a universal family (see \cite[Prop.8.8.2, p.233]{BL04}), and in particular, so does $V^\p$. Therefore, it suffices to prove that $M(U^\p/V^\p)$ surjects onto $W$. By \cite{Do92}, there is a surjection $\pi_1(V) \twoheadrightarrow W$ whose kernel $N:=\ker(\pi_1(V) \twoheadrightarrow W)$ is the image of $\pi_1(\ti{V})$ under the natural homomorphism.  The group $M:= \lambda_*(\pi_1(V^\p))= \ker(\pi_1(V) \ra G)$ determines the Galois cover $V^\p \ra V$.  The situation is summarized in the following commutative diagram
$$\xymatrix{
\pi_1(\ti{U}^\p) \ar[r] \ar[d] & \pi_1(\ti{U}) \ar[d] \\
\pi_1(V^\p) \ar[r]^-{\lambda_*} & \pi_1(V) \ar@{->>}[d] \\
&  W, 
}$$
where the cokernel of the vertical map on the left is the monodromy group $M(U^\p/V^\p)$. To prove surjectivity of $M(U^\p/V^\p) \ra W$ induced by the homomorphisms in the above diagram, it suffices to show that $M$ surjects onto $W$. By considering the indices of $M$ and $N$ in $\pi_1(V)$ we see immediately that $N$ is not contained in $M$. Indeed, the index of $M$ is the order of $G$ and is given by a standard formula. However, it suffices to note that $H = \mrm{GL}(5, \mb{Z}/3)/(\pm I)$ injects as a subgroup of $G$ and $|H| = {\frac 1 2} (3^5 - 1) (3^5 - 3) (3^5 - 3^2) (3^5 - 3^3) (3^5 - 3^4)$, which is already much larger than 51840. Therefore, the image $\bar{N}$ of $N$ in the quotient $\pi_1(V)/M \simeq G$ is non-trivial. Since $N$ is a normal subgroup of $\pi_1(V)$, $\bar{N}$ is normal in the quotient group $G$.  Since $G$ is a simple group (cf. \cite[Thm.1, p.12]{Di48}) and $\bar{N}$ is non-trivial, we must have $\bar{N}=G$. Thus, we have $MN = \pi_1(V)$, which shows that the quotient map $M \ra \pi_1(V)/N \simeq W$ is surjective.  
 
By Lemma \ref{lem6}, there exists a smooth connected curve section $S$ of $V$ whose topological fundamental group surjects onto that of $V$. Let $S^\p :=\lambda^{-1}(S)$ and let $\mc{R}^\p$ be the pull-back of $U^\p$ to $S^\p$. The covering $\mu \colon S^\p \ra S$ is determined by the inverse image $Q$ of $\lambda_*(\pi_1(V^\p))$ in $\pi_1(S)$ under the homomorphism $\pi_1(S) \ra \pi_1(V)$ induced by the inclusion (in particular, $Q=\mu_*\pi_1(S^\p)$). Therefore, since $\lambda_*$ is injective, from the following commutative diagram 
$$\xymatrix{
\pi_1(V^\p) \ar@{^{(}->}[r]^-{\lambda_*} & \pi_1(V) \\
\pi_1(S^\p) \ar[r]^-{\mu_*} \ar[u] & \pi_1(S) \ar@{->>}[u] 
}$$
we see that the vertical map on the left is surjective. Using the following commutative diagram 
$$\xymatrix{
\pi_1(\ti{U}^\p) \ar[r] &  \pi_1(V^\p)\\
\pi_1(\ti{\mc{R}}^\p) \ar[r] \ar[u] &  \pi_1(S^\p) \ar[u]
}$$ 
and surjectivity of $\pi_1(S^\p) \twoheadrightarrow \pi_1(V^\p)$   we conclude that the monodromy group of $\mc{R}^\p \ra S^\p$ still surjects onto $W$. 

To summarize some of the constructions so far, we now have over the smooth curve $S'$, a degree 27 finite \'etale map $\mc{R}^\p \to S^\p$ and a family $\mc{P}^\p \to S^\p$ of ppav's (with level 3 structure), such that for each point $s^\p \in S^\p$, each one of the 27 points in the fiber of $\mc{R}^\p \to S^\p$ over $s^\p$ (equivalently, in the fiber in $\ms{R}_6$ over the image point of $V \subset \ms{A}_5$) corresponds to exactly one of 27 Abel-Prym curves in the abelian variety $\mc{P}^\p_{s^\p}$.  Each Abel-Prym curve is uniquely determined up to translations (since the automorphism group of the ppav $\mc{P}^\p_{s^\p}$ is $\{\pm 1\}$, and an Abel-Prym curve always has some translates that are invariant under $-1$).  Since translation in the abelian variety defines an obvious algebraic equivalence between a 1-cycle and any of its translates, maintaining a careful distinction between an Abel-Prym curve $C_i$ and the collection of all its translates, is not the essential point in our arguments; we are going to suppress the distinction, for simplicity.  What is critical is that the mapping $\mc{R}^\p \to S^\p$ still has monodromy group $W$.

Additionally, for simplicity of notation, we shall now write $\mc{R}$ and $S$ for $\mc{R}^\p$ and $S^\p$, respectively, in the sequel.  Thus, with this updated notation, let $\mc{P}$ be the universal family of 5-dimensional ppav's (with level 3 structure) over $S$.  The fiber over $s \in S$ will be denoted by $\mc{P}_s$.   

Let $\ms{V}$ be the local system on $S$ whose fiber over a point $s \in S$ is the $27$-dimensional vector space spanned by the Abel-Prym curves in $\mc{P}_s$. Given a contractible open set $U \subset S$ in the analytic topology we fix a trivialization 
$$\ms{V}(U) \simeq \bigoplus_{i=1}^{27}\mb{Q}\mc{C}_i,$$ 
and let $C_{i,s}$ denote the fiber of $\mc{C}_i$ over $s \in U$. Thus, the fiber of $\ms{V}$ over $s\in U$ is $\ms{V}_s = \bigoplus_{i=1}^{27}\mb{Q}C_{i,s}$. To simplify the notation, for $q=(q_1, \ldots, q_{27}) \in \mb{Q}^{27}$ define 
$$Z_s(q):=\sum_{i=1}^{27} q_i C_{i,s}.$$
Since $S$ has a base for the analytic topology consisting of contractible sets, the assignment 
$$U \mapsto \ms{T}(U):=\big\{Z_s(q) \,|\, Z_s(q) \sim_{\mrm{alg}} 0 \text{ for all } s \in U\big\}$$
defines a presheaf on $S$ in the analytic topology. In the sequel, we let $\ms{T}$ denote the associated sheaf on $S$.  

\begin{lem4} \label{lem4}
The subsheaf $\ms{T}$ of $\ms{V}$ is a local subsystem of $\ms{V}$. 
\end{lem4}
\begin{proof} Given a contractible open set $U \subset S$ in the analytic topology, a trivialization $\ms{V}(U) \simeq \bigoplus_{i=1}^{27}\mb{Q}\mc{C}_i$ and $q = (q_1, \ldots, q_{27}) \in \mb{Q}^{27}$ as above, we claim that the locus
$$\{s\in U \,|\, Z_s(q) \sim_{\mrm{alg}} 0 \} $$
is either all of $U$ or a countable set of points in $U$. 
Assuming this claim, let us prove the lemma. Fix $s \in S$ for the remainder of the proof and consider the monodromy representation $$\rho \colon \pi_1(S,s) \ra \mrm{Aut}(\ms{V}_s).$$
It is well-known that there is a natural bijection between isomorphism classes of local systems and monodromy representations up to conjugation, see \cite[Cor.3.10, p.71]{V07}. Under this bijection subrepresentations correspond to local subsystems. Therefore, to show that $\ms{T}$ is a local subsystem of $\ms{V}$, it suffices to prove that for each $[\gamma] \in \pi_1(S,s)$ the automorphism $\rho([\gamma])$ of $\ms{V}_s$ maps $\ms{T}_s$ to itself. The image of  $\gamma \colon [0,1] \ra S$ can be covered by finitely many contractible open sets $U_1, \ldots, U_n$ for the analytic topology. Assume $s \in U_1$ and $U_i \cap U_{i+1} \ne \emptyset$ for $i=1,\ldots, n-1$. Given an element $Z_s(q) \in \ms{T}_s$, by the above claim we know that $Z_t(q) \in \ms{T}_t$ for all $t\in U_1$. Applying the same argument to $U_2,\ldots, U_n$,  we conclude that the parallel transport of $Z_s(q)$ stays in the fibers of $\ms{T}$. Therefore, $\rho([\gamma])$ maps $\ms{T}_s$ to itself.

In the remainder of the proof we show the above claim. Let $\ti{\mc{R}} \ra S$ be the Galois closure of the cover $\mc{R} \ra S$.  Let $\ti{\ms{V}}$ be the pull-back of $\ms{V}$ to $\ti{\mc{R}}$ under the composition $\ti{\mc{R}} \ra \mc{R} \ra S$. Since the claim is about algebraic equivalences in the fibers of $\ms{V}$, it suffices to show the claim for $\ti{\ms{V}}$ globally over $\ti{\mc{R}}$. The local system $\ti{\ms{V}}$ trivializes over $\ti{\mc{R}}$ and we may write 
$$\ti{\ms{V}} \simeq \bigoplus_{i=1}^{27}\mb{Q}\mc{C}_i.$$
Let $\ti{\mc{P}}:= \mc{P} \times_S \ti{\mc{R}} \ra \ti{\mc{R}}$ be the pull-back of the universal family. 
By clearing the denominators and rearranging the terms we may rewrite $Z_s(q) \sim_{\mrm{alg}} 0$ as
\begin{equation} \label{eq1}
\sum_{i \in I} a_iC_{i,s} \sim_{\mrm{alg}} \sum_{j \in J} b_jC_{j,s},
\end{equation} where $a_i, b_j$ are non-negative integers, $I \cup J = \{1, \ldots, 27\}$ and $I\cap J = \emptyset$.  Furthermore, $(\ref{eq1})$ holds if and only if there exists an effective cycle $E_s$ on $\ti{\mc{P}}_s$ such that the cycles $E_s+\sum_{i \in I} a_iC_{i,s}$ and $E_s+\sum_{j \in J} b_jC_{j,s}$ lie in the same connected component of the Chow variety of $\ti{\mc{P}}_s$, see \cite[4.1.3, p.122]{K96}. 
In what follows we may and shall assume that the effective cycle $E_s$ is effectively algebraically equivalent to a sufficiently high multiple of the cycle $\Xi^4_s$, where $\Xi_s$ is the theta divisor on $\ti{\mc{P}}_s$ (see \cite[4.1.2, p.121]{K96} for the definition of effective algebraic equivalence). More precisely, we may arrange, by taking another finite cover if necessary, that there is a family over $\ti{\mc{R}}$ of theta divisors $\Xi_s \subset \ti{\mc{P}}_s$, and we take $E_s$ to be the intersection of four general elements of the linear system $|n\Xi_s|$ for $n \gg 0$. To emphasize the role of $n$ we write $E_{n,s}$ for such $E_s$. This assumption is needed for the following two reasons. First, to ensure that $E_{n,s}$ extends to a cycle over all of $\ti{\mc{R}}$, i.e., that there is an effective cycle $\mc{E}_n$ flat over $\ti{\mc{R}}$ whose fiber over $s \in \ti{\mc{R}}$ is $E_{n,s}$. Second, to ensure that is suffices to consider at most countably many effective relative cycles $\mc{E}_n$ (indexed by $n$) in order to exhibit algebraic equivalence (provided it holds) of $\sum_{i \in I} a_iC_{i,t}$ and  $\sum_{j \in J} b_jC_{j,t}$ on $\ti{\mc{P}}_t$ by adding on the fiber of $\mc{E}_{n}$ over $t$. 



Consider the relative Chow variety of $\ti{\mc{P}}$ over $\ti{\mc{R}}$ that parametrizes effective cycles on $\ti{\mc{P}}$ which lie over $0$-dimensional subschemes of $\ti{\mc{R}}$, see \cite[Def.3.1.1, p.41 and Def.3.20, p.51]{K96}. Let $E_{n,t}$ denote the fiber of $\mc{E}_n$ over $t\in \ti{\mc{R}}$. It is easily seen that the cycles $E_{n,t} + \sum_{i \in I} a_iC_{i,t}$ and $E_{n,t} + \sum_{j \in J} b_jC_{j,t}$ are  algebraically equivalent on $\ti{\mc{P}}$ (but not necessarily on $\ti{\mc{P}}_t$), and therefore, lie in the same connected component $\mc{X}$ of the relative Chow variety of $\ti{\mc{P}}$ over $\ti{\mc{R}}$ . Applying Proposition \ref{pro1} to the natural morphism $\mc{X} \ra \ti{\mc{R}}$ we conclude that the set of $t \in \ti{\mc{R}}$ such that $E_{n,t} +\sum_{i \in I} a_iC_{i,t}$ and $E_{n,t} +\sum_{j \in J} b_jC_{j,t}$ lie in the same connected component of the Chow variety of $\ti{\mc{P}}_t$ is Zariski closed in $\ti{\mc{R}}$. 

If $\sum_{i \in I} a_iC_{i,t}$ and $\sum_{j \in J} b_jC_{j,t}$ are algebraically equivalent on $\ti{\mc{P}}_t$, there exists $m \ge n$ and a relative cycle $\mc{E}_m$ (constructed as above) such that $E_{m,t} + \sum_{i \in I} a_iC_{i,t}$ and $E_{m,t} + \sum_{j \in J} b_jC_{j,t}$ lie in the same component of the Chow variety of $\ti{\mc{P}}_t$. Furthermore, since the choice of the four effective divisors in $|m\Xi_t|$ that define $E_{m,t}$ is unrestricted (as long as they intersect properly), it suffices to consider at most countably many such cycles $\mc{E}_m$ in order to exhibit algebraic equivalence of $\sum_{i \in I} a_iC_{i,t}$ and $\sum_{j \in J} b_jC_{j,t}$ on $\ti{\mc{P}}_t$ (for any $t$). Therefore, we conclude that the locus of $t \in \ti{\mc{R}}$ such $\sum_{i \in I} a_iC_{i,t}$ and $\sum_{j \in J} b_jC_{j,t}$ are algebraically equivalent on $\ti{\mc{P}}_t$ is either all of $\ti{\mc{R}}$ or a countable union of points in $\ti{\mc{R}}$.
\end{proof}

\begin{thm1} \label{thm1}
The $\mb{Q}$-vector space of algebraic equivalences of the $27$ Abel-Prym curves on a generic principally polarized abelian 5-fold has dimension $20$.
\end{thm1}
\begin{proof}
By Lemma \ref{lem3} the fibers of $\ms{T}$ have dimension $\ge 20$. By definition of $\ms{T}$, for each $s \in S$ the 1-dimensional sub-representation $\ms{L}_s$ of $\ms{V}_s$ spanned by the  Beauville degree $0$ graded piece of an Abel-Prym curve intersects $\ms{T}_s$ trivially. By Lemma \ref{lem1}, the composition $\ms{T}_s \hookrightarrow \ms{V}_s \ra \ms{V}_s/\ms{L}_s$ is not surjective for generic $s \in S$. Therefore,  for generic $s\in S$ we have $\dim \ms{T}_s\le 25$. By Lemma \ref{lem4}, $\ms{T}$ is invariant under the monodromy of the Galois closure of $\mc{R} \ra S$, which by our construction surjects onto the Weyl group $W$. By Lemma \ref{lem5} we conclude that the fibers of $\ms{T}$ have dimension $20$. 
\end{proof}

\begin{cor1}
On a generic ppav $P$ of dimension $5$ with Abel-Prym curves $\ti{C}_1,\ldots, \ti{C}_{27}$ we have:
(1) $\ms{T}(P, \ti{C}_i) \ne \ms{T}(P, \ti{C}_j)$ if $i \ne j$; 
(2) $[\ti{C}_i]_{(2)} \ne 0$ in $A^*(P)$ for $1\le i \le 27$.
\end{cor1}
\begin{proof} (1) For each $i \in \{1, \ldots, 27\}$, the bi-graded piece of $\ms{T}(P, \ti{C}_i)$ of codimension 4 and Beauville degree 2 is spanned by $[\ti{C}_i]_{(2)}$. If $\ms{T}(P, \ti{C}_i) = \ms{T}(P, \ti{C}_j)$ then $[\ti{C}_i]_{(2)}$ and $[\ti{C}_j]_{(2)}$ are proportional, which is easily seen to contradict the structure of the vector space of relations among Abel-Prym curves from the proof of  Theorem \ref{thm1}. 

(2) We may check that the class $[\ti{C}_i]_{(2)}$ does not belong to the 20-dimensional $\mb{Q}$-vector space of algebraic equivalences among $\ti{C}_1, \ldots, \ti{C}_{27}$ from Theorem \ref{thm1}.
\end{proof}

\section*{Acknowledgements}
The first named author expresses his gratitude to David Swinarski for help with {\tt Magma}. Both authors are thankful to Mitchell Rothstein for a helpful suggestion. The first named author is grateful to the second named author and to Roy Smith for sharing a question by Emanuele Raviolo that is closely related to the subject of the present article, which was then in progress.

\end{document}